\documentclass[12pt]{amsart}

\usepackage{amsmath}
\usepackage{amssymb}
\usepackage{amscd}
\usepackage[mathscr]{eucal}

\newtheorem{theorem}{Theorem}[section]
\newtheorem{lemma}[theorem]{Lemma}
\newtheorem{corollary}[theorem]{Corollary}
\newtheorem{proposition}[theorem]{Proposition}

\theoremstyle{definition}

\newtheorem{example}[theorem]{Example}

\theoremstyle{remark}

\newtheorem{notation}[theorem]{Notation}

\numberwithin{equation}{section}

\newcommand{\N}{\mathbb{N}}

\newcommand{\real}{\mathbb{R}}

\newcommand{\s}{\mathbb{S}}
\newcommand{\B}{\mathbb{B}}
\newcommand{\J}{\mathbb{R}_+}
\newcommand{\Cl}{{\rm Cl}}
\newcommand{\Int}{{\rm Int}}

\newcommand{\asdim}{{\rm asdim}}
\newcommand{\cat}{{\rm CAT(0)}}

\newcommand{\diam}{{\rm diam\,}}

\usepackage[dvips]{graphicx,color}

\begin{document}

\title%[Preliminary manuscript]
{Asymptotic dimension of proper CAT(0) spaces which are homeomorphic to the plane}

\author{Naotsugu Chinen$^{*}$ \and Tetsuya Hosaka}

\address{Okinawa National College of Technology,
Nago-shi Okinawa 905-2192, Japan}
 
\email{chinen@okinawa-ct.ac.jp}

\address{Department of Mathematics, Faculty of Education, 
Utsunomiya University, Utsunomiya, 321-8505, Japan}

\email{hosaka@cc.utsunomiya-u.ac.jp}

 

\thanks{ 
$^{*}$Partly supported by the Grant-in-Aid for Scientific Research (C),
Japan Society for the Promotion of Science, No. 19540108.
}

\begin{abstract} 
In this paper, we investigate 
a proper CAT(0) space $(X,d)$ 
which is homeomorphic to $\real^2$ and 
we show that the asymptotic dimension $\asdim (X,d)$ is 
equal to $2$.
\end{abstract}
 
\maketitle

\section{Introduction and preliminaries}
\medskip

In this paper, we study asymptotic dimension of proper CAT(0) spaces 
which are homeomorphic to $\real^2$.

A metric space $(X,d)$ is {\it proper}
if all closed, bounded sets in $(X,d)$ are compact.
We say that a metric space $(X,d)$ is a {\it geodesic space} if 
for any $x,y \in X$, 
there exists an isometric embedding $\xi:[0,d(x,y)] \rightarrow X$ such that 
$\xi(0)=x$ and $\xi(d(x,y))=y$ (such a $\xi$ is called a {\it geodesic}).

Let $(X,d)$ be a geodesic space and 
let $T$ be a geodesic triangle in $X$.
A {\it comparison triangle} for $T$ is 
a geodesic triangle $\overline{T}$ in the Euclidean plane $\real^2$
with same edge lengths as $T$.
Choose two points $x$ and $y$ in $T$. 
Let $\bar{x}$ and $\bar{y}$ denote 
the corresponding points in $\overline{T}$.
Then the inequality $$d(x,y) \le d_{\real^2}(\bar{x},\bar{y})$$ 
is called the {\it CAT(0)-inequality}, 
where $d_{\real^2}$ is the usual metric on $\real^2$.
A geodesic space $X$ is called a {\it CAT(0) space} 
if the CAT(0)-inequality holds
for all geodesic triangles $T$ and for all choices of two points $x$ and 
$y$ in $T$.
Details of CAT(0) spaces are found in \cite{BH}.

In Section~2, 
we first investigate 
proper CAT(0) spaces which are homeomorphic to $\real^2$ 
and we show the following.

\begin{proposition}
Let $(X,d)$ be a proper $\cat$ space which is homeomorphic to $\real^2$.
Then,
$S(x,r)$ is homeomorphic to $\s^1$ for all $x \in X$ and all $r > 0$.
Hence the boundary $\partial X$ is homeomorphic to a circle $\s^1$.
\end{proposition}

Let $(X,d)$ be a metric space and 
let $\mathscr{U}$ be a family of subsets of $(X,d)$.
The family $\mathscr{U}$ is said to be {\it uniformly bounded}
if there exists a positive number $K$
such that
$\diam U \leq K$ 
for all $U \in \mathscr{U}$.
The family $\mathscr{U}$ is said to be {\it $r$-disjoint}
if $d(U,U')>r$ for any $U,U'\in \mathscr{U}$ with $U \neq U'$.

The {\it asymptotic dimension} of a metric space 
$(X,d)$ does not exceed $n$ and we write $\asdim (X,d) \leq n$ 
if for every $r > 0$ there exist uniformly bounded,
$r$-disjoint families $\mathscr{U}^0,\mathscr{U}^1,\ldots,\mathscr{U}^n$
of subsets of $X$ such that $\bigcup_{k=0}^n\mathscr{U}^k$ covers $X$.
The {\it asymptotic dimension} of a metric space $(X,d)$ is equal to $n$
and we write $\asdim (X,d) = n$ if $\asdim (X,d) \leq n$ 
and $\asdim (X,d) \not\leq n-1$.

Asymptotic dimension of groups 
relates to the Novikov conjecture and 
there are some interesting recent research on asymptotic dimension
(cf.\ \cite{BD}, \cite{DJ}, \cite{DKU}, \cite{Gr}, \cite{Yu}).
In \cite{Gr}, Gromov remarks that 
word hyperbolic groups have finite asymptotic dimension 
and Roe gave detail of the proof in \cite{R}.
Asymptotic dimension of CAT(0) groups and CAT(0) spaces is unknown 
in general.

The purpose of this paper is to prove the following theorem.

\begin{theorem}\label{MainThm}
If $(X,d)$ is a proper $\cat$ space 
which is homeomorphic to $\real^2$, 
then $\asdim(X,d) = 2$.
\end{theorem}

We note that the proper CAT(0) space $(X,d)$ in this theorem 
need not have an action of some group.
We give an example in Section~4.

\medskip
\section{Proper $\cat$ spaces which are homeomorphic to $\real^2$}
\medskip

We first give notation used in this paper.

\begin{notation}
Let the set of all natural number, real number and $[0,\infty)$
denote by $\N$, $\real$ and $\J$, respectively.
Set $\J^n = \real^{n-1} \times \J$,
$\B^n = \{x \in \real^{n} : \sum_{i=1}^{n} x_i^2 \leq 1\}$ and 
$\s^n = \{x \in \real^{n+1} : \sum_{i=1}^{n+1} x_i^2 =1\}$.
Let $Y$ be a subspace of a metric space $(X,d)$.
The interior and the closure of $Y$ in a space $X$
will be denoted by $\Int_{X}Y$ and $\Cl_{X}Y$, respectively.
Also set $B(x,r) = \{y \in X : d(x,y) \leq r\}$ and
$S(x,r) = \{y \in X : d(x,y) = r\}$.
We denote the geodesic from $x$ to $y$ in a CAT(0) space $(X,d)$ 
by $[x,y]$ (cf. \cite[Proposition II 1.4]{BH}).
Set $[x,y) = [x,y] \setminus \{y\}$,
$(x,y] = [x,y] \setminus \{x\}$ and
$(x,y) = [x,y] \setminus \{x, y\}$.
\end{notation}

The following lemma is known.

\medskip

\begin{lemma}\label{lemma:comp}
Let $(X,d)$ be a proper $\cat$ space,
$r > 0$ and $x_0 \in X$.
Then, the following are satisfied:
\begin{enumerate}
\item 
$B(x_0,r)$ is a convex set.
\item 
$x_0 \not\in [x,y] \subset B(x_0,r)$ and $(x,y) \subset B(x_0,r) \setminus S(x_0,r)$
for any $x, y \in S(x_0,r)$ with $d(x,y) < 2r$.
\item$($cf.\cite[Lemma II 5.8 and Proposition II 5.12]{BH}$)$
If $X$ is a manifold,
for each $x \in X \setminus \{x_0\}$
there exists a geodesic line $\xi : \real \rightarrow X$ such that
$\xi(0) = x_0$ and $\xi(d(x_0,x)) = x$.
\end{enumerate}
\end{lemma}

\medskip

We investigate a proper $\cat$ space which is homeomorphic to $\real^2$.

\medskip

\begin{notation}\label{notation:ell}
Let $(X,d), r, x_0, x,$ and $y$ be as in Lemma \ref{lemma:comp}(2).
Suppose that $X$ is homeomorphic to $\real^2$.
By Lemma \ref{lemma:comp},
there exist two geodesic rays $\xi_{x_0,x}, \xi_{x_0,y} : \J \rightarrow X$ such that
$\xi_{x_0,x}(0) = \xi_{x_0,y}(0) =x_0$,
$\xi_{x_0,x}(r) = x$ and $\xi_{x_0,y}(r) = y$.
By Lemma \ref{lemma:comp},
$\xi_{x_0,x}([r,\infty))\cup [x,y] \cup \xi_{x_0,y}([r,\infty))$ is homeomorphic to $\real$.
Since $X$ is homeomorphic to $\real^2$,
by Sch$\ddot{\rm o}$nflies Theorem,
there exists the component $C$ of $X \setminus \xi_{x_0,x}([r,\infty))\cup [x,y] \cup \xi_{x_0,y}([r,\infty))$
such that
$x_0 \not\in C$.
Set $\ell(x,y) = S(x_0,r) \cap \Cl_X C$.
\end{notation}

\medskip

We show some lemmas.

\medskip

\begin{lemma}\label{lemma:b}
Let $(X,d)$ be a proper $\cat$ space which is homeomorphic to $\real^2$.
Then,
$S(x,r)$ is a continuum for all $x \in X$ and all $r > 0$.
\end{lemma}

\begin{proof}
Let $x_0 \in X$ and $r > 0$.
Since $B(x_0,r)$ is a convex set,
by duality (cf. \cite{S}),
\[\bar H_0 (X \setminus B(x_0,r)) \cong \check H^1(B(x_0,r)) = 0,\]
thus, $X \setminus \Int_X B(x_0,r) = \Cl_X (X \setminus B(x_0,r))$ is connected.
Since there exists a deformation retraction of $X \setminus \Int_X B(x_0,r)$ onto $S(x_0,r)$,
$S(x_0,r)$ is connected.
\end{proof}

\medskip

\begin{lemma}\label{lemma:ab}
Let $(X,d)$ be a proper $\cat$ space which is homeomorphic to $\real^2$,
$r > 0$, $x_0 \in X$, $x, y \in S(x_0,r)$ with $0 < d(x,y) < 2r$
and $z \in \ell(x,y)$.
Then 
\begin{enumerate}
\item
$\ell(x,y)$ is a continuum,

\item
$[x,y] \cap [x_0,z] \neq \emptyset$ and

\item
$d(x,z) \leq d(x,y)$.
\end{enumerate}
\end{lemma}

\begin{proof}
(1)
By Notation \ref{notation:ell},
there exists the component $D$ of $X \setminus R$
such that
$C \subset D$,
where $R = \xi_{x_0,x}(\J) \cup \xi_{x_0,y}(\J)$.
Since $X$ is homeomorphic to $\real^2$,
by Sch$\ddot{\rm o}$nflies Theorem,
$\Cl_X D$ is homeomorphic to $\J^2$.
Let $D'$ be a copy of $D$.
Define an equivalent relation $\sim$ in $D \cup D'$ as follows:
for $a \in D$ and $a' \in D'$,
$a \sim a'$ if and only if $a=a'$, $a \in R$ and $a' \in R'$.
Set $B = B(x_0,r)\cap \Cl_X D$, 
$\widetilde{D} = (D \cup D')/\sim$
and 
$\widetilde{B} = (B \cup B')/\sim$.
Then, 
there exists a deformation retraction $\Cl_X(D \setminus B)$ onto $\ell(x,y)$,
$\widetilde{D}$ is homeomorphic to $\real^2$
and
$\widetilde{B}$ is a contractible compact set.
By the same method as in the proof of Lemma \ref{lemma:b},
we can show that $\ell(x,y) \cup (\ell(x,y))'/\sim$ is connected.
Since there exists the natural surjective map from $\ell(x,y) \cup (\ell(x,y))'/\sim$ onto $\ell(x,y)$,
$\ell(x,y)$ is connected.

(2)
We may assume that $z \not\in \{x,y\}$.
By Notation \ref{notation:ell},
there exists the component $C$ of $X \setminus \xi_{x_0,x}([r,\infty))\cup [x,y] \cup \xi_{x_0,y}([r,\infty))$
such that
$x_0 \not\in C$ and $z \in C$.
Thus, $\xi_{x_0,x}([r,\infty))\cup [x,y] \cup \xi_{x_0,y}([r,\infty))$ separates between $x_0$ and $z$ in $X$.
Since $[x_0,z] \subset B(x_0,r)$ is an arc connecting between $x_0$ and $z$ in $X$,
$[x,y] \cap [x_0,z] \neq \emptyset$.

(3)
On the contrary, suppose $d(x,z) > d(x,y)$.
By (2), 
there exists $z' \in [x,y] \cap [x_0,z]$.
Since $z' \in [x,y]$,
\[ d(x,z') + d(z',y) = d(x,y) < d(x,z) \leq d(x,z') + d(z',z),\]
thus, $d(z',y) < d'(z',z)$.
Then, 
\[r = d(x_0,y) \leq d(x_0,z') + d(z',y) < d(x_0,z')+ d(z',z) =d(x_0,z) = r,\]
a contradiction.
\end{proof}

\medskip

\begin{lemma}\label{lemma:conn}
Let $(X,d)$ be a proper $\cat$ space which is homeomorphic to $\real^2$,
$r,t > 0$, $x_0 \in X$ and  $y_0 \in S(x_0,r)$.
Then $S(x_0,r) \cap B(y_0,t)$ is connected.
\end{lemma}

\begin{proof}
Set $N = S(x_0,r) \cap B(y_0,t)$.
If $t \geq 2r$,
$S(x_0,r) \subset B(y_0,t)$.
By Lemma \ref{lemma:b},
$N$ is connected.
We may assume that $t < 2r$.
Take $x \in N$.
Since $d(y_0,x) \leq t < 2r$,
by Lemma \ref{lemma:ab},
we have
\[\ell(y_0,x) \subset S(x_0,r) \cap B(y_0,d(y_0,x)) \subset S(x_0,r) \cap B(y_0,t) = N.\]
Therefore,
by Lemma \ref{lemma:ab},
$N = \bigcup \{\ell(y_0,x) : x \in N\}$ is connected,
which proves the lemma.
\end{proof}

\medskip

We obtain the following proposition from lemmas above.

\begin{proposition}\label{proposition:R^2}
Let $(X,d)$ be a proper $\cat$ space which is homeomorphic to $\real^2$.
Then,
$S(x,r)$ is homeomorphic to $\s^1$ for all $x \in X$ and all $r > 0$.
\end{proposition}

\begin{proof}
By Lemma \ref{lemma:b} and \cite[Theorem 11.21]{W},
it suffices to show the following:
\begin{enumerate}
\item
$S(x_0,r) \setminus \{y_0,y_1\}$ is non-connected
for any $y_0,y_1 \in S(x_0,r)$ with $y_0 \neq y_1$.

\item
$S(x_0,r) \setminus \{y_0\}$ is connected
for each $y_0 \in S(x_0,r)$.
\end{enumerate}

We take two points $y_0, y_1 \in S(x_0,r)$ with $y_0 \neq y_1$.
By Lemma \ref{lemma:comp}, 
there exist geodesic rays $\xi_{x_0,y_0}, \xi_{x_0,y_1} : \J \rightarrow X$
such that
$\xi_{x_0,y_0}(0)=\xi_{x_0,y_1}(0) = x_0$,
$\xi_{x_0,y_0}(r) = y_0$ and $\xi_{x_0,y_1}(r) =y_1$.
By Sch$\ddot{\rm o}$nflies Theorem,
there exist closed sets $Z_0,Z_1$ of $X$ 
such that
$Z_i$ is homeomorphic to $\J^2$ for $i =0,1$, 
$X = Z_0 \cup Z_1$ and $Z_0 \cap Z_1 \subset \xi_{x_0,y_0}(\J) \cup \xi_{x_0,y_1}(\J)$
is homeomorphic to $\real$.
Since $S(x_0,r) \cap \Int_X Z_i \neq \emptyset$ for $i = 0,1$,
$S(x_0,r)\setminus \{y_0,y_1\}$ is non-connected,
which proves (1).

Let $x,y \in S(x_0,r) \setminus \{y_0\}$ with $x \neq y$.
By Lemma \ref{lemma:comp}, 
there exist geodesic rays $\xi_{x_0,x}, \xi_{x_0,y} : \J \rightarrow X$
such that
$\xi_{x_0,x}(0)= \xi_{x_0,y}(0)= x_0$,
$\xi_{x_0,x}(r) = x$ and $\xi_{x_0,y}(r) = y$.
Set $R = \xi_{x_0,x}(\J) \cup \xi_{x_0,y}(\J)$.
By \cite[Proposition 1.4(1), p.160]{BH},
there exists $z \in [x_0,x)$ such that
$\xi_{x_0,x}(\J) \cap \xi_{x_0,y}(\J) = [x_0,z]$.
By Sch$\ddot{\rm o}$nflies Theorem,
there exists the component $C$ of $X \setminus R$ 
such that
$y_0 \not\in C$ and
$E_{x,y} = \Cl_X C$ is homeomorphic to $\J^2$.
Set $L_{x,y}= E_{x,y}\cap S(x_0,r)$.
We see $B(x_0,d(x_0,z)) \subset E_{x,y}$ or $B(x_0,d(x_0,z)) \cap E_{x,y} = \{z\}$.
Suppose that $B(x_0,d(x_0,z)) \subset E_{x,y}$.
We note that 
$L_{x,y}$, $B(x_0,d(x_0,z))$ and $\{z\}$
are deformation retracts of $\Cl_X(E_{x,y} \setminus B(x_0,r))$, $E_{x,y} \cap B(x_0,r)$
and 
$B(x_0,d(x_0,z))$, respectively.
Thus, $\{z\}$ is a deformation retract of $E_{x,y} \cap B(x_0,r)$.
By the same method as in the proof of Lemma \ref{lemma:ab} (1),
we can show that 
$L_{x,y}\cup (L_{x,y})'/\sim$ is a deformation retract of 
$\Cl_X(E_{x,y} \setminus B(x_0,r))\cup(\Cl_X(E_{x,y} \setminus B(x_0,r)))'/\sim$ 
and 
$\{z\}$ is a deformation retract of 
$(E_{x,y} \cap B(x_0,r))\cup(E_{x,y} \cap B(x_0,r))'/\sim$,
thus,
$L_{x,y}$ is connected.
Suppose that $B(x_0,d(x_0,z)) \cap E_{x,y} = \{z\}$.
Since $\{z\}$ is a deformation retract of $E_{x,y} \cap B(x_0,r)$,
by the same method above,
we can show that $L_{x,y}$ is connected.

Fix $y_0' \in S(x_0,r) \setminus \{y_0\}$.
Since $S(x_0,r)\setminus \{y_0\} = \bigcup \{L_{x,y_0'} : x \in S(x_0,r) \setminus \{y_0,y_0'\}\}$,
it is connected,
which proves (2).
\end{proof}

\medskip

\begin{corollary}
If $(X,d)$ is a proper $\cat$ space which is homeomorphic to $\real^2$, 
then the boundary $\partial X$ of $X$ is homeomorphic to $\s^1$.
\end{corollary}

\medskip

We show the following lemma which is used in the proof of the main theorem.

\medskip

\begin{lemma}\label{lemma:arc}
Let $(X,d)$ be a proper $\cat$ space which is homeomorphic to $\real^2$,
$x_0 \in X$,  $r,t > 0$ with $2t < r$ and
$x, x' \in S(x_0,r)$ with $3t \leq d(x,x') < 2r$.
Then there exist $y_0, \dots, y_{3n-1} \in S(x_0,r)$ and $m \in \N$ with $0 < 3m < 3n-1$ 
such that
$y_0 = x$, $y_{3m} = x'$,
$t \leq \diam \ell(y_i,y_{i+1}) \leq 2t$,
$\{y_0,\dots, y_{3n-1}\} \cap \ell(y_i,y_{i+1}) = \{y_i,y_{i+1}\}$
for each $i = 0,\dots, 3n-1$,
$S(x_0,r) = \ell(y_0,y_{1})\cup \cdots \cup \ell(y_{3n-1},y_{3n})$
and
$\ell(x,x') = \ell(y_0,y_{1})\cup \cdots \cup \ell(y_{3m-1},y_{3m})$,
where $y_{3n} =y_0$.
\end{lemma}

\begin{proof}
Set $z_0 = y_0 = x$.
By Proposition \ref{proposition:R^2},
$S(x_0,r)$ is homeomorphic to $\s^1$.
Since $S(x_0,r) \not\subset B(z_0,t)$,
by Lemma \ref{lemma:conn},
$S(x_0,r) \cap B(z_0,t)$, $\ell(x,x')$ and $\ell(x,x') \cap B(z_0,t)$ are arcs.
Let $z_1$ be the end point of $\ell(x,x') \cap B(z_0,t)$ with $z_0 \neq z_1$.
By Lemma \ref{lemma:conn},
we have $\ell(z_0,z_1)= \ell(x,x') \cap B(z_0,t)$.
By Lemma \ref{lemma:conn},
$S(x_0,r) \cap B(z_1,t)$ is an arc.
Since $z_0$ and $z_1$ are the end points of $\ell(z_0,z_1)$ with $d(z_0,z_1) = t$,
there exists the end point $z_2$ of $S(x_0,r) \cap B(z_1,t)$
such that
$\diam \ell(z_1,z_2) = t$,
$\ell(z_0,z_1)\cap \ell(z_1,z_2) = \{z_1\}$
and $\ell(z_0,z_1)\cup \ell(z_1,z_2) = \ell(x,x') \cap B(z_1,t)$.
Thus, by induction,
we can take $z_2,\dots,z_{p+1} \in S(x_0,r)$
and 
an arc $\ell(z_{i-1},z_i)$ in $\ell(x,x')$ with the end points $\{z_{i-1},z_i\}$
such that
$z_{p+1} = x'$,
$\ell(z_{i-1},z_i)\cap \ell(z_i,z_{i+1}) = \{z_i\}$ 
for each $i = 1,\dots,p$,
$\ell(z_{i-1},z_i)\cup \ell(z_i,z_{i+1}) = \ell(x,x') \cap B(z_{i},t)$,
for each $i = 1,\dots,p$,
$\ell(x,x') = \bigcup_{i=1}^{p+1} \ell(z_{i-1},z_{i})$,
$\diam \ell(z_{i-1},z_{i}) = t$ for any $i=1,\dots,p$ and
$\diam \ell(z_p,z_{p+1}) \leq t$.
Let $k \in \N$ and $\delta = 0,1,2$ such that $p = 3k + \delta$.
Set $m =k$ and $y_{3m} = z_{p+1}$.
If $\delta = 0$,
set $y_i = z_i$ for each $i = 1,\dots,3m-1$.
If $\delta = 1$,
set $y_i = z_{i}$ for each $i = 1,\dots,3m-2$ and $y_{3m-1} = z_{p-1}$.
If $\delta = 2$,
set $y_i = z_{i}$ for each $i = 1,\dots,3(m-1)$,
$y_{3m-2} = z_{p-3}$ and $y_{3m-1} = z_{p-1}$.
Similarly, we have $y_{3m+1}, \dots, y_{3n-1} \in \Cl_X(S(x_0,r) \setminus \ell(x,x'))$,
which proves the lemma.
\end{proof}

\medskip

\medskip
\section{Asymptotic dimension of 
proper $\cat$ spaces which are homeomorphic to $\real^2$}
\medskip

First we show the following.

\begin{lemma}\label{lemma:2over}
Let $(X,d)$ be a proper \cat~ space which is homeomorphic to $\real^2$.
Then, $\asdim(X,d) \geq 2$.
\end{lemma}

\begin{proof}
On the contrary, suppose that $\asdim (X,d) \leq 1$.
Let $r > 0$. 
There exist uniformly bounded,
$3r$-disjoint families $\mathscr{U}^0,\mathscr{U}^1$
of subsets of $X$ such that $\mathscr{U}^0 \cup \mathscr{U}^1$ covers $X$.
Since $X$ is homeomorphic to $\real^2$,
there exist uniformly bounded,
$r$-disjoint families $\mathscr{V}^0,\mathscr{V}^1$
of subsets of $X$ 
satisfying the following:
\begin{enumerate}
\item
$\mathscr{V}^0 \cup \mathscr{V}^1$ covers $X$.
\item
Every $V \in \mathscr{V}^0 \cup \mathscr{V}^1$ is a compact topological 2-manifold with boundary.
\end{enumerate}
Let $\varepsilon > 0$ with $\varepsilon < r/2$,
let $V \in \mathscr{V}^i$ and
let $M$ and $M'$ be two components of $V$ with $d(M,M') = d(M, V \setminus M) < \varepsilon$.
Then there exists a disk $A$ in $X$
such that
$M \cup A \cup M'$ is connected,
$V \cup A$ is a compact topological 2-manifold with boundary
and
$d(V \cup A,V') > r-\varepsilon$ whenever $V' \in \mathscr{V}^i$ with $V \neq V'$.
Thus, we may assume that
\begin{enumerate}
\item[(3)]
$d(M,M') \geq \varepsilon$
for each $V \in \mathscr{V}^0 \cup \mathscr{V}^1$ and each two components $M,M'$ of $V$.
\end{enumerate}
Since $\mathscr{V}^0 \cup \mathscr{V}^1$ is uniformly bounded,
there exists 
$r \leq s = \sup \{\diam C : C$ is a component of $V \in \mathscr{V}^0 \cup \mathscr{V}^1\} < \infty$.
Thus, we may assume that there exists a component $C_0$ of $V_0 \in \mathscr{V}^0$
such that $s- \varepsilon < \diam C_0 \leq s$.
We have $c_{0},c_{1} \in C_0$ such that $d(c_{0},c_{1}) = \diam C_0$.
By Lemma \ref{lemma:comp},
there exists a geodesic line $\xi : \real \rightarrow X$
such that
$\xi (0) = c_{0}$ and $\xi (\diam C_0) = c_{1}$.
Since $C_0 \cap \xi(\real) \subset \xi([0,\diam C_0])$,
there exists the component $N$ of $\partial C_0$ containing $c_{0},c_{1}$
which is contained in the closure of the unbounded component of $X \setminus C_0$.

We note that $N \subset \bigcup \{ V_1 : V_1 \in \mathscr{V}^1\}$
is homeomorphic to $\s^1$.
Since $\mathscr{V}^1$ is $r$-disjoint,
there exists a component $C_1$ of $V_1 \in \mathscr{V}^1$ 
such that $N \subset C_1$.
Then, there exist $t_0, t_1 \in \real$ with $t_0 < 0 < \diam C_0 < t_1$
such that $\xi(t_0), \xi(t_1) \in C_1 \cap \xi(\real) \subset \xi([t_0,t_1])$.
By the similar argument above,
we can show there exist a component $N'$ of $\partial C_1$ containing $\xi(t_0), \xi(t_1)$
and
$V_2 \in \mathscr{V}^0$ containing $N'$.
If $V_0 = V_2$,
by (3),
$d(c_0,\xi(t_0)) \geq \varepsilon$ and $d(c_1,\xi(t_1)) \geq \varepsilon$, i.e.,
$d(\xi(t_0),\xi(t_1)) > \diam C_0 + 2\varepsilon > s$, which contradicts the definition of $s$.
If $V_0 \neq V_2$,
$d(V_0,V_2) > r$.
Thus,
$d(c_0,\xi(t_0)) > r$ and $d(c_1,\xi(t_1)) > r$, i.e.,
$d(\xi(t_0),\xi(t_1)) > \diam C_0 + 2r > s$, which contradicts the definition of $s$.
\end{proof}

\medskip

We prove the main theorem.

\medskip

\begin{theorem}\label{theorem:cat}
Let $(X,d)$ be a proper $\cat$ space which is homeomorphic to $\real^2$.
Then, $\asdim(X,d) = 2$.
\end{theorem}

\begin{proof}
By Lemma \ref{lemma:2over}
it suffices to show that $\asdim(X,d) \leq 2$.

Let $r > 0$.
Fix $x_0 \in X$ and $k \in \N$ with $k \geq 6$.
By Lemma \ref{lemma:arc},
there exist $y_{0,0}, \dots, y_{0,3n(0)-1} \in S(x_0,kr)$
such that
$2r \leq \diam \ell(y_{0,i},y_{0,i+1}) \leq 16r$,
$\{y_{0,0},\dots, y_{0,3n(0)-1}\} \cap \ell(y_{0,i},y_{0,i+1}) = \{y_{0,i},y_{0,i+1}\}$
for each $i = 0,\dots, 3n(0)-1$ and
$S(x_0,kr) = \ell(y_{0,0},y_{0,1})\cup \cdots \cup \ell(y_{0,3n(0)-1},y_{0,3n(0)})$,
where $y_{0,3n(0)} =y_{0,0}$.
See Figure \ref{theorem:cat}.1.
Set $\mathscr{V}_{0,\delta} = \{\ell(y_{0,3i+\delta},y_{0,3i+1+\delta}) : i = 0,\dots,n(0)-1\}$
for each $\delta = 0,1,2$.

%%%%%%%%%% picture %%%%%%%%
\begin{center}
\setlength{\unitlength}{3pt}
\begin{picture}(120,62)
\put(60,10){\circle*{1}}\put(60,7){$x_0$}

\put(0,60){\footnotesize Im $\xi _{0,i}$}
\qbezier[160](60,10)(50,40)(0,60)
\put(20,50.8){\circle*{1}}\put(15,53){\footnotesize $y'_{0,i}$}
\put(38.5,38.7){\circle*{1}}\put(35.5,40.5){\footnotesize $y_{0,i}$}

\put(115,60){\footnotesize Im $\xi _{0,i+1}$}
\qbezier[160](60,10)(70,40)(120,60)
\put(100,50.8){\circle*{1}}\put(97,53){\footnotesize $y'_{0,i+1}$}
\put(81.5,38.7){\circle*{1}}\put(82,40){\footnotesize $y_{0,i+1}$}

\put(0,27){\footnotesize $S(x_0,kr)$}
\qbezier(60,40)(30,40)(0,30)
\qbezier(60,40)(90,40)(120,30)

\put(0,42){\footnotesize $S(x_0,(k+1)r)$}
\qbezier(60,55)(30,55)(0,45)
\qbezier(60,55)(90,55)(120,45)

\qbezier[160](60,10)(60,10)(60,55)
\put(60,55){\circle*{1}}\put(60,57){$y'_{0,i,0}$}
\put(60,40){\circle*{1}}\put(60,42){$y_{0,i,0}$}
\put(54,27){\footnotesize $[x_0,y'_{0,i,0}]$}

\put(50,1){Figure \ref{theorem:cat}.1}
\end{picture}
\end{center}
%%%%%%%%%% picture %%%%%%%%

For every $i = 0,\dots, 3n(0)-1$ there exists a geodesic ray 
$\xi_{0,i} : \J \rightarrow X$ such that
$\xi_{0,i}(0) = x_0$ and $\xi_{0,i}(kr) = y_{0,i}$.
Set $y_{0,i}' = \xi_{0,i}((k+1)r)$ for each $i = 0,\dots, 3n(0)-1$.
We note $2r \leq d(y_{0,i}, y_{0,i+1}) < d(y_{0,i}', y_{0,i+1}') \leq 18r$
for each $i = 0,\dots, 3n(0)-1$.

Let $i \in \{0,\dots,3n(0)-1\}$.

Suppose that 
$d(y_{0,i}', y_{0,i+1}') < 12r$.
We can take $y_{0,i,0} \in \ell(y_{0,i}, y_{0,i+1})$ and $y_{0,i,0}' \in \ell(y_{0,i}', y_{0,i+1}')$
such that 
$r \leq d(y_{0,i},y_{0,i,0})= d(y_{0,i+1},y_{0,i,0})=d(y_{0,i},y_{0,i+1})/2 < 6r$ and
$\{y_{0,i,0}\} = [x_0,y_{0,i,0}']\cap S(x_0,kr)$.
We note $r < d(y_{0,i}',y_{0,i,0}'), d(y_{0,i+1}',y_{0,i,0}') < 8r$.

Suppose that 
$12r \leq d(y_{0,i}', y_{0,i+1}')$.
We note that $10r \leq d(y_{0,i}, y_{0,i+1})$.
There exist $z_{0,i,0}, z_{0,i,1} \in \ell(y_{0,i},y_{0,i+1})$ and $z_{0,i,0}', z_{0,i,1}' \in \ell(y_{0,i}',y_{0,i+1}')$
such that
$d(y_{0,i},z_{0,i,0}) = d(y_{0,i+1},z_{0,i,1}) = r$ and
$\{z_{0,i,j}\} = [x_0,z_{0,i,j}']\cap S(x_0,kr)$ for $j = 0,1$.
We note that $d(y_{0,i}',z_{0,i,0}'), d(y_{0,i+1}',z_{0,i,1}') \leq 3r$
and 
$6r \leq d(z_{0,i,0}', z_{0,i,1}')$.
By Lemma \ref{lemma:arc},
there exist $y_{0,i,1}', \dots, y_{0,i,3k_{0,i}-1}' \in \ell(y_{0,i,0}',y_{0,i,3k_{0,i}}')$
such that
$2r \leq d(y_{0,i,j}',y_{0,i,j+1}') \leq 4r$
and
$\ell(y_{0,i,j}',y_{0,i,j+1}')\cap \{y_{0,i,0}', \dots, y_{0,i,k_{0,i}}'\} = \{y_{0,i,j}',y_{0,i,j+1}'\}$
for each $j = 0,\dots,3k_{0,i}-1$,
where $y_{0,i,0}' = z_{0,i,0}'$ and $y_{0,i,3k_{0,i}}' = z_{0,i,1}'$.
See Figure \ref{theorem:cat}.2.

%%%%%%%%%% picture %%%%%%%%
\begin{center}
\setlength{\unitlength}{3pt}
\begin{picture}(120,62)
\put(60,10){\circle*{1}}\put(60,7){$x_0$}

\put(0,60){\footnotesize Im $\xi _{0,i}$}
\qbezier[160](60,10)(50,40)(0,60)
\put(20,50.8){\circle*{1}}\put(15,53){\footnotesize $y'_{0,i}$}
\put(38.5,38.7){\circle*{1}}\put(35.5,40.5){\footnotesize $y_{0,i}$}

\put(115,60){\footnotesize Im $\xi _{0,i+1}$}
\qbezier[160](60,10)(70,40)(120,60)
\put(100,50.8){\circle*{1}}\put(101.5,51.5){\footnotesize $y'_{0,i+1}$}
\put(81.5,38.7){\circle*{1}}\put(82,40){\footnotesize $y_{0,i+1}$}

\put(0,27){\footnotesize $S(x_0,kr)$}
\qbezier(60,40)(30,40)(0,30)
\qbezier(60,40)(90,40)(120,30)

\put(0,42){\footnotesize $S(x_0,(k+1)r)$}
\qbezier(60,55)(30,55)(0,45)
\qbezier(60,55)(90,55)(120,45)

\put(30,25){\footnotesize $[x_0,y'_{0,i,0}]$}\put(42,29){\vector(1,1){6}}
\qbezier[120](60,10)(50,40)(30,52)
\put(30,52.5){\circle*{1}}\put(25,54.5){\footnotesize $y'_{0,i,0}$}
\put(45,39.5){\circle*{1}}\put(43.5,41.5){\footnotesize $z_{0,i,0}$}

\put(77,25){\footnotesize $[x_0,y'_{0,i,3k_{0,i}}]$}\put(78,29){\vector(-1,1){6}}
\qbezier[120](60,10)(70,40)(90,52)
\put(90,52.5){\circle*{1}}\put(90,55){\footnotesize $y'_{0,i,3k_{0,i}}$}
\put(75,39.5){\circle*{1}}\put(70,41.5){\footnotesize $z_{0,i,1}$}

\put(38,53.6){\circle*{1}}\put(33.3,56.5){\footnotesize $y'_{0,i,1}$}
\put(45,54.6){\circle*{1}}\put(42.3,57.5){\footnotesize $y'_{0,i,2}$}
\put(60,57){$\cdots$}
\put(80,54){\circle*{1}}\put(74,56.5){\footnotesize $y'_{0,i,3k_{0,i}-1}$}

\put(50,1){Figure \ref{theorem:cat}.2}
\end{picture}
\end{center}
%%%%%%%%%% picture %%%%%%%%

Set $Y_1 = \{y_{0,i,j}' : 0\leq i \leq 3n(0)-1 \mbox{~and~} j = 0,\dots,3k_{0,i}\}$
and
$n(1) \in \N$ with $3n(1) -1 = |Y_1|$.
We can renumber $Y_1 = \{y_{1,i} : i = 0,\dots,3n(1)-1\}$
such that
\begin{align*}
\{y' \in S(x_0,(k+1)r) : y \in \bigcup \mathscr{V}_{0,1}\cap \bigcup \mathscr{V}_{0,2}\}
\subset \bigcup\mathscr{V}_{1,0},\\
\{y' \in S(x_0,(k+1)r) : y \in \bigcup \mathscr{V}_{0,0}\cap \bigcup \mathscr{V}_{0,2}\}
\subset \bigcup\mathscr{V}_{1,1},\\
\{y' \in S(x_0,(k+1)r) : y \in \bigcup \mathscr{V}_{0,0}\cap \bigcup \mathscr{V}_{0,1}\}
\subset \bigcup\mathscr{V}_{1,2},
\end{align*}
and
$\ell(y_{1,i},y_{1,i+1}) \cap Y_1 = \{y_{1,i},y_{1,i+1}\}$
for each $i = 0,\dots,3n(1)-1$,
where 
we let 
$y_{1,3n(1)} = y_{1,0}$ and
$\mathscr{V}_{1,\delta} = \{\ell(y_{1,3i+\delta},y_{1,3i+1+\delta}) : i = 0,\dots,n(1)-1\}$
for each $\delta = 0,1,2$.
We note that $2r \leq \diam V \leq 16r$ for all $\delta = 0,1,2$ and all $V \in \mathscr{V}_{1,\delta}$.

By induction,
for every $m \in \N$ with $m \geq 2$
there exists
$Y_m = \{y_{m,i} : i = 0,\dots,3n(m)-1\} \subset S(x_0,(k+m)r)$
such that
\begin{align*}
\{y' \in S(x_0,(k+m)r) : y \in \bigcup \mathscr{V}_{m-1,1}\cap \bigcup \mathscr{V}_{m-1,2}\}
\subset \bigcup\mathscr{V}_{m,0},\\
\{y' \in S(x_0,(k+m)r) : y \in \bigcup \mathscr{V}_{m-1,0}\cap \bigcup \mathscr{V}_{m-1,2}\}
\subset \bigcup\mathscr{V}_{m,1},\\
\{y' \in S(x_0,(k+m)r) : y \in \bigcup \mathscr{V}_{m-1,0}\cap \bigcup \mathscr{V}_{m-1,1}\}
\subset \bigcup\mathscr{V}_{m,2},
\end{align*}
$\ell(y_{m,i},y_{m,i+1}) \cap Y_m = \{y_{m,i},y_{m,i+1}\}$
for each $i = 0,\dots,3n(m)-1$
and
$2r \leq \diam V \leq 16r$ for all $\delta = 0,1,2$ and all $V \in \mathscr{V}_{m,\delta}$,
where
we let $\mathscr{V}_{m,\delta} = \{\ell(y_{m,3i+\delta},y_{m,3i+1+\delta}) : i = 0,\dots,n(m)-1\}$
for each $\delta = 0,1,2$.

Set 
$\overline{V} = \{x \in B(x_0,(k+m+1)r) \setminus \Int_X B(x_0,(k+m)r) : 
[x_0,x] \cap V \neq \emptyset\}$
for each $V \in \mathscr{V}_{m,\delta}$,
$\overline{\mathscr{V}_{m,\delta}} = \{\overline{V} : V \in \mathscr{V}_{m,\delta}\}$
and
$\mathscr{W}_\delta = \left\{ W : W \mbox{~is a component of 
$\bigcup_{m=0}^\infty \overline{\mathscr{V}_{m,\delta}}$}\right\}$
for each $\delta = 0,1,2$.
By construction,
we have the following:
\begin{enumerate}
\item 
For $V \in \mathscr{V}_{m,\delta}$, 
$\diam \overline{V} \cap S(x_0,(k+m+1)r) < 12r$
if and only if
$\overline{V} \cap \bigcup \overline{\mathscr{V}_{m+1,\delta}} = \emptyset$.

\item
Let $\mathscr{V}_{m+1}(V) =
\{U \in \mathscr{V}_{m+1,\delta} : 
\overline{V} \cap U \neq \emptyset\}$
for each $V \in \mathscr{V}_{m,\delta}$.
Then, $U \subset \overline{V}$ 
for $U \in \mathscr{V}_{m+1}(V)$.

\item
We have $\mathscr{V}_{m+2}(U) = \emptyset$ for each $V \in \mathscr{V}_{m,\delta}$ and each $U \in \mathscr{V}_{m+1}(V)$
because $\diam U < 12r$ by construction.
\end{enumerate}
For every $\delta = 0,1,2$ and every $W \in \mathscr{W}_\delta$
we have 
\begin{align*}
\diam W & \leq \sup\{\diam V : V \in \mathscr{V}_{m,\delta}
\mbox{~for $m \geq 0$ and $\delta = 0,1,2 $\}} + 4r \\
&\leq 16r + 4r = 20r.
\end{align*}

Let $V_i,V_j \in \mathscr{V}_{m,\delta}$ with $V_i \neq V_j$.
We show that $d(V_i,V_j) \geq r$.
On the contrary, suppose that 
$d(x,y) < r$ for some $x \in V_i$ and some $y \in V_j$.
By Lemma \ref{lemma:conn},
let $\ell(x,y)$ denote the arc in $S(x_0,(k+m)r) \cap B(x,d(x,y))$ with the end points $\{x,y\}$.
By construction,
we have $i = 0,\dots,n(m)-1$ such that $\ell(y_{m,i},y_{m,i+1}) \subsetneq \ell(x,y)$.
However, 
$r \leq \diam \ell(y_{m,i},y_{m,i+1}) \leq \diam \ell(x,y) = d(x,y) < r$,
a contradiction.

Let $\overline{V_i}, \overline{V_j} \in \overline{\mathscr{V}_{m,\delta}}$ 
with $\overline{V_i} \neq \overline{V_j}$.
We show that $d(\overline{V_i},\overline{V_j}) \geq r$.
Let $x' \in \overline{V_i}$ and $y' \in \overline{V_j}$.
Set $\{x\} = [x_0,x'] \cap V_i$ and $\{y\} = [x_0,y'] \cap V_j$.
By above, $r \leq d(V_i,V_j) \leq d(x,y)$.
Let $T$ be the geodesic triangle consisting of
three points $x_0,x',y'$, 
let $\overline{T}$ be a comparison triangle for $T$ in $\real^2$ and
let $\overline{x_0}, \overline{x}, \overline{y}, \overline{x'},$ and $\overline{y'}$ 
denote the corresponding points in $\overline{T}$.
Since $X$ is a $\cat$ space,
we have 
\[r \leq d(x,y) \leq d_{\real^2}(\overline{x}, \overline{y}) 
\leq d_{\real^2}(\overline{x'}, \overline{y'}) = d(x',y'),\]
thus, $d(\overline{V_i},\overline{V_j}) \geq r$.

Let $\overline{V_i}\in \overline{\mathscr{V}_{m,\delta}}$ and 
$\overline{V_j} \in \overline{\mathscr{V}_{m+1,\delta}}$ 
with $\overline{V_i} \cap \overline{V_j} = \emptyset$.
Set 
$W_j = \{[x_0,x] \cap S(x_0,(k+m)) : x \in \overline{V_j}\}$.
By the definition of $y_{m,i,j}'$'s,
similarly, we can show $d(V_i,W_j) \geq r$.
Since $X$ is a $\cat$ space,
by the same method
we can obtain that $d(\overline{V_i}, \overline{V_j}) \geq r$.
By (1), (2) and (3),
we have $d(W,W') \geq r$
for any $W,W' \in \mathscr{W}_\delta$ with $W \neq W'$.

Let
$\mathscr{U}_0 = \left\{ U : U \mbox{~is a component of 
$B(x_0,kr) \cup \bigcup \mathscr{W}_{0}$}\right\}$
and
$\mathscr{U}_\delta = \mathscr{W}_{\delta}$
for $\delta =1,2$.
By above,  $\mathscr{U}_0 \cup \mathscr{U}_1 \cup \mathscr{U}_2$ is 
a uniformly bounded cover of $(X,d)$ and
$d(U,U') \geq r$
for any $U,U' \in \mathscr{U}_{\delta}$ with $U \neq U'$,
which proves the theorem.
\end{proof}

\medskip

\medskip
\section{Application}
\medskip

As an application of Theorem~\ref{theorem:cat}, 
we obtain the following corollary.

\medskip

\begin{corollary}
Let $(W,S)$ be a Coxeter system.
If the boundary $\partial\Sigma(W,S)$ of $\Sigma(W,S)$ is homeomorphic to $\s^1$, 
then $\asdim W = 2$.
\end{corollary}

\begin{proof}
Let $(W,S)$ be a Coxeter system 
whose boundary $\partial\Sigma(W,S)$ is homeomorphic to $\s^1$.
Then the Coxeter group $W$ is a virtual Poincar\'{e} duality group, 
and for some $\tilde{S}\subset S$, 
$$W=W_{\tilde{S}}\times W_{S\setminus \tilde{S}},$$
where the nerve $N(W_{\tilde{S}},\tilde{S})$ is homeomorphic to $\s^1$ 
and $W_{S\setminus \tilde{S}}$ is finite (\cite{D3}, cf.\ \cite{H}).
Then the Davis complex $\Sigma(W,S)$ splits as 
$$\Sigma(W,S)=\Sigma(W_{\tilde{S}},\tilde{S})
\times\Sigma(W_{S\setminus \tilde{S}},S\setminus \tilde{S}).$$
Here 
$\Sigma(W_{\tilde{S}},\tilde{S})$ is homeomorphic to $\real^2$ and 
$\Sigma(W_{S\setminus \tilde{S}},S\setminus \tilde{S})$ is bounded.
By Theorem~\ref{theorem:cat}, 
we obtain that $\asdim \Sigma(W,S) = 2$.
Hence $\asdim W = 2$.
\end{proof}

\medskip

In general, it is known that 
every Coxeter group has finite asymptotic dimension 
(\cite{DJ}, cf.\ \cite{DS}).

\medskip

\begin{example}\label{example:mn}
Let 
$m \in \N$ 
and
let $D_m \subset \real^2$ be a regular $m$-polygon
with a metric $d_m = d_{\real^2}|_{D_m}$ and edges $e_1,\dots,e_m$
such that $\diam e_i = 1$ for each $i = 1,\dots, n$.
We consider a noncompact cell 2-complex $(\Sigma,d)$ with a triangulation $\mathscr{T}$ as follows.
\begin{enumerate}
\item 
For every $\sigma \in \mathscr{T} \setminus \mathscr{T}^{(1)}$
there exist $m(\sigma) \in \N$ and a simplicial isometry $f_{\sigma}$ from $(D_{m(\sigma)},d_{m(\sigma)})$ onto $(|\sigma|,d|_{|\sigma|})$.
\item
$|\{\sigma  \in \mathscr{T} \setminus \mathscr{T}^{(1)} : \tau < \sigma\}| = 2$
for each $\tau \in \mathscr{T}^{(1)} \setminus \mathscr{T}^{(0)}$.
\item
$r(v) = \sum \{\pi - 2\pi/m(\sigma) : v < \sigma \in \mathscr{T} \setminus \mathscr{T}^{(1)}\} \geq 2\pi$
for each $v \in \mathscr{T}^{(0)}$.
\item
for any $x,y \in \Sigma$
\begin{align*}
d(x,y) &= \min\{ \sum_{j=1}^k d_{m(\sigma)}(f_{\sigma_j}^{-1}(x_{j-1}),f_{\sigma_j}^{-1}(x_j)) : \\
&\hspace{2cm} x= x_0 \in |\sigma_1|, x_j \in |\sigma_j| \cap |\sigma_{j+1}| (1 \leq j < k), 
y = x_{k} \in |\sigma_{k}| \}.\\
\end{align*}
\end{enumerate}

By \cite{D2}, every $(\Sigma,d)$ above is
a $\cat$ space which is homeomorphic to $\real^2$, hence
we obtain that $\asdim (\Sigma,d) = 2$ from Theorem \ref{theorem:cat}.
Here we note that 
$(\Sigma,d)$ need not have an action of some group and 
$(\Sigma,d)$ is neither a euclidean nor a hyperbolic plane in general.
\end{example}

\end{document}